\documentclass[journal]{IEEEtran}
\usepackage{diagbox}
\usepackage{amsfonts}
\usepackage{cite}
\usepackage{graphicx}
\usepackage{amsmath}
\usepackage{amsthm}
\usepackage{amssymb}
\usepackage{bm}
\ifCLASSINFOpdf
\else
\fi
\begin{document}
%\title{Elastic Load Shaping for Power Distribution Grid with Local Area Packetized Power Networks}
%\title{Electricity Pricing beyond Accumulative Perspective}
%\title{Algebraic Insights on Elements of Power System}
\title{Algebraic Analysis of\\ Electric Load and Electricity Pricing}

%\title{Electricity Pricing over $L_2$ Space}
%State Key Laboratory of Advanced Optical \\Communication Systems and Networks
%\author{\IEEEauthorblockN{Jinghuan Ma}
%\IEEEauthorblockA{School of Electronics Engineering \\and Computer Science\\
%Peking University\\
%Beijing, China 100871\\
%Email: mjhdtc@pku.edu.cn}
%\and
%\IEEEauthorblockN{Lingyang Song}
%\IEEEauthorblockA{School of Electronics Engineering \\and Computer Science\\
%Peking University\\
%Beijing, China 100871\\
%Email: lingyang.song@pku.edu.cn}}

\author{
Jinghuan~Ma, Jie~Gu, and Zhijian~Jin

\thanks{
J.~Ma, J.~Gu and Z.~Jin are with the School of Electronic Information and Electrical Engineering,
Shanghai~Jiao~Tong~University, Shanghai, China, 200240 (email: mjhdtc@sjtu.edu.cn).
}
}


%Telephone: (800) 555--1212, Fax: (888) 555--1212}


\maketitle

\begin{abstract}
This paper presents a perspective of functional analysis to analyze electric load and electricity pricing in the $L_2$ space and its isomorphic vector space, which also yields the general algebraic model of load and pricing associated with time course. It clarifies law of interaction between payment and load associated with time course and necessity of pricing to sufficiently convey information on how load results in supply cost. It describes classic paradigm of modeling in a general algebraic form and formally proves the ineffectiveness of classic integral-based pricing in modeling. It consequently introduces to:
describe electricity supply cost by a mapping defined on an isomorphic space of the original space of load constituted by orthonormal basis that quantifies the dynamism of load, named the space of dynamism; a pricing model defined on the space of dynamism that sufficiently conveys the information;
a further derived computational model of pricing based on the Fourier series;
simple examples to demonstrate use of the proposed pricing and its effectiveness in distinguishing loads and reflecting supply cost associated with time course, which theoretically completes the incentive to reshape loads.


%Simulation results confirm that
%the proposed algorithm significantly improves consumer¡¯s cost
%efficiency compared with the conventional payoff maximization
%algorithm. It also shows that a higher service fee will considerably
%decrease the cost efficiency.
\end{abstract}


% \fi
\IEEEpeerreviewmaketitle

%\vspace{-1em}
\section{Introduction}
In a general sense, electric load, cost of electricity supply and electricity pricing, as elements describing the supply-demand interaction of power system~\cite{PSE-2018,MOE-2002}, can be classified into two classes. Electric load and cost of electricity supply, as physical phenomena, can be interpreted as \emph{hard} natural interactions; while electricity pricing, as a force that does not conducted by matter, can be compared to \emph{soft} artificial impact. But the soft impact does not really come from nowhere, as to complete feedback loop of
the supply-demand interaction, electricity pricing must convey to the demand side the cost of supply and payment must be determined based on the load.

Since the beginning of power industry, people have been, partially, pushing the supply-demand interaction toward economic optimality, system stability and sustainability, by improving the
\emph{soft} impact. Specifically, the improvement refers to improving the capability in reflecting true facts of the \emph{hard} phenomena, i.e., how electric load results in cost of electricity supply. The flat rate~\cite{ET-1951} is the earliest pricing and has a single function to simply reflect a rough and averaged cost of supply. Soon, the multiple tariffs\footnote{The concept refers to different flat rates subject to the purpose of supply~\cite{ET-1951}.}~\cite{ET-1951} and the two-part tariff~\cite{ET-1951,TT-1941} were proposed and adopted to address the fact that loads in different cases~(also interpreted as loads in different forms, different loads) resulting in different supply costs.
It is noteworthy that the concept of marginal cost~\cite{SO-1948,RD-1949}, which has been widely accepted and continuously used till today despite the changes in specific form, has been adopted by the two-part tariff at that time. This shows that people have attempted, at such an early stage, to improve the effectiveness of pricing by improving its capability in distinguishing loads and restoring relationships between loads and supply costs. The terms \emph{discrimination} and \emph{discriminatory} used by William Vickrey in his early work~\cite{SO-1948} are equivalent expressions of \emph{distinguishing}; and the work suggests that discrimination plays an important role in resource allocation. Though it has not further explained the capability of discrimination mathematically.

Model is built to describe, explain and manipulate phenomena. Things that cannot be captured cannot be manipulated; things that cannot be distinguished cannot be manipulated separately.
In this paper, the capability of distinguishing elements of same kind and relationships between elements of different kinds is presented as the capability of modeling in a mathematical sense of generality~\cite{EO-1989}.

For electric load model, to distinguish elements of same kind refers to abstractly describing the relationships among electric loads, so as to distinguish an element from one another in accordance with that a load distinguishes itself from one another in the real world. For supply cost model, to distinguish relationships between elements of different kinds refers to abstractly describing the causality of load and supply cost by mapping, so as to distinguish a cost-load causal relation from one another in accordance with the real causality of electric load and supply cost that different loads and loads matched by different combinations of generation units results in non-negligible difference in supply cost.
Electricity pricing has to address both: to distinguish loads and charge them based on the supply cost caused, for fairness and preciseness~(direction and strength of impact).
The capability of modeling reaches a complete level once the model can completely distinguish abstract elements and relationships, strictly and completely corresponding to the phenomena in reality. That is, a bijective mapping can be defined between the model and the reality.

A key factor that restrains paradigm of modeling is the capability in observation and computation. A simple case is that it requires at least two points to fix a linear function, while it requires more sufficient observations to validate the linear approximation.
At the early stage, it was impossible to automatically detect phenomenon, record data, and process data with high accuracy, let alone with high detecting frequency. Amount and diversity of information observed on features and relationships of the phenomena was extremely limited, which could not afford to build a model of high information capacity.

The classic paradigm of modeling load, cost characteristic, and supply-demand relationship~\cite{ET-1951,TT-1941,SO-1948,RD-1949,PG-2013} had been customized since then. That is, to outline a time lasting phenomenon by a model independent of time. To completely describe demand of electricity, one has to model it as a function of time, which is an element that belongs to a linear space of real-valued square integrable functions~($L_2$ space)~\cite{EO-1989}. The $L_2$ space is a Hilbert space of infinite dimensions, which is far more complex than a field of real numbers. But it was even difficult to record an hourly load profile for every load in the 1930s to 1950s. People could only use real-numbered or integer-numbered labels to distinguish loads~\cite{SO-1948}, such as service pattern, size of daily demand, use pattern of certain marginal cost. As supply-demand activity is time-continuous, the relationship between supply cost and load lies in time course of any length, instead of time course of certain length or being independent of time. However, due to the extremely limited capability in observation and computation, the cost-load relationship could only be described by time-independent curves of two dimensions~\cite{ET-1951,SO-1948,PG-2013}, of which the amount and diversity of information, had been dramatically compressed compared to the original relationship that has to be described in the $L_2$ space. The approximation of an $F-P$ curve~\cite{SO-1948,PG-2013} used tens or hundreds of data pairs of fuel costs and outputs, which means the cost characteristic of a generator, that originally lies in time course of any length, had be outlined by limited information contained in the tens or hundreds of observations that were not associated with time course at all.

People have being driven to strengthen the capability of modeling by insufficiency of the older models compared to observations of the time that have been continuously improved.
The most straightforward and most widely adopted approach has been to decompose a long time course into multiple shorter cycles.
Time-of-use tariff~\cite{TO-2012} has been proposed since the 1950s and has evolved from time-of-day tariff~\cite{ET-1951} to peak-load pricing~\cite{PL-1957,PL-1960}, and to spot pricing~(also known as real-time pricing)~\cite{OP-1982,EW-2003}. The length of the cycle has been shortened from a day~\cite{ET-1951}, to a couple of hours~\cite{PL-1957,PL-1960}, an hour~\cite{OP-1982}, and five minutes~\cite{EW-2003}. In each cycle, models of the classic paradigm that are independent of time have been elaborated to describe marginal costs w.r.t. different causes~(generation cost, transmission loss, congestion price, etc.) and different locations, etc~\cite{OP-1982,EW-2003}. Other factors associated with pricing, such as welfare,
utility and payoff~\cite{OP-1982,AD-2010,RR-2011}, are also time-independent in each cycle.

Today, advancements in observation, computation and communication have being creating huge amount of data that contains abundant information on features and relationships of power system, with high accuracy and high sampling frequency. Volume and quality of data are keeping increasing. It is intuitive to shorten the length of cycle to a more real-time level.
But it has been reported both in China and the U.S. that real-time pricing in the wholesale market cannot effectively reflect the real supply cost~\cite{RR-2011,DT-2018,RO-2020}.
Objections to marginal-cost pricing that follows the classic paradigm of modeling, raised by William Vickery~\cite{SO-1948}, can date back to 1948.

The cause cannot simply be explained as inevitable gap between reality and model.

The gap is believed to be caused by the classic paradigm of modeling.
By decomposing a long time course into multiple shorter cycles, we break the continuity of the operation in time domain and break the integrity
of information on the physical properties associated with undivided time course of electricity consumption of any length.
Furthermore, by defining variables and relationships independent of time in each cycle, the classic paradigm of modeling actually compress the original phenomena
as a time course into a snapshot, like a frame that outlines a time-lasting motion. A vast majority of the original information is undoubtedly lost.
Thirdly, it is worth a rethink whether we really utilize the information contained in the fine-grained data. We use it as input in each cycle to a predetermined static model of relationship that does not require extra information.
%The output of the model does not rely on real characteristic of the system, but associated with what has been used to predetermine the model of relationship.
We have not used it to improve the modeling capability in capturing features and relationships associated with time course.
For example, in each cycle of locational marginal pricing, what represents the relationship between generation cost and load is the coefficients of marginal generation cost that are predetermined according to unit characteristics of generators.
%The unit characteristics are usually $F-P$ curves approximated by a limited number of data obtained in former steady operations, which have nothing to do with new information and can hardly restore the real states.
The output only reflects an assumed outcome of the unit characteristics based on experience of the steady test circumstances.
In this paper, classic paradigm of modeling is described in a general algebraic form; and its insufficiency is formally proved.

Judgement yielded by the former discussion is that mismatch between information capacity of model~(capability of paradigm of modeling) and information capacity of data~(capability in observation and computation) has reached a condition that appeals for breakthrough in paradigm of modeling.

With algebraic expressions addressing generality, we presents a perspective of functional analysis to
1)~analyze electric load and electricity pricing in the $L_2$ space and its isomorphic vector space, which also yields the general model of load and pricing associated with time course;
2)~clarify law of interaction between payment and load associated with time course and necessity of pricing to sufficiently convey information on how load results in cost;
3)~describe classic paradigm of modeling in a general algebraic form and formally prove the ineffectiveness of classic integral-based pricing in modeling.
We consequently introduce
1)~to describe supply cost by a mapping defined on an isomorphic space of the original space of load constituted by orthonormal basis that quantifies the dynamism of load, named the space of dynamism;
2)~a pricing model defined on the space of dynamism that sufficiently conveys the information;
3)~a computational model of load and pricing derived from the general forms in the space of dynamism based on the Fourier series;
4)~simple examples to demonstrate use of the proposed pricing and its effectiveness in distinguishing loads and reflecting supply cost associated with time course, which theoretically completes the incentive to reshape load curves.

The remainder of this paper is organized as follows.
General model of electric load and analysis of load modeling is presented in Section~\ref{CGM}.
General model of electricity pricing and analysis of classic modeling is presented in Section~\ref{A2}.
In Section~\ref{PRPS}, models of load, supply cost and pricing defined on the space of dynamism as isomorphic vector space of the original space of load are introduced; the further derived computational model of load and pricing based on the Fourier series is introduced; examples of the proposed pricing are presented and analyzed.
Conclusion is drawn in Section~\ref{clnc}.

\section{General Model of Electric Load and\\ Analysis of Load Modeling}\label{CGM}

In this section, appearing natures of electric load are translated into abstract expression in the form of algebra without distortion of physical meanings.
With the natures abstracted as sufficient conditions, the space of load functions, as the abstract form of the set of loads, are identified as a linear space of real-valued square integrable functions~($L_2$ space), such that more properties of the $L_2$ space and its elements can be applied to the space of load functions, which provides significant inspirations on modeling, analyzing and impacting electric load.
With the general model as reference, algebraic analysis of classic paradigm of load modeling is presented.

\vspace{-0.5em}
\subsection{Abstraction of Electric Load}

%Let $E$ denote the energy accumulated on a time span $\Delta t$.
%Classically, the rate of energy consumption, also known as power, is defined as the derivative of energy consumed with respect to time:
%$dE \over dt$.
%However, in this work, we define a kind of average power as
%\begin{equation}
%p(t,{\Delta}t)={E\big|_{ [t-{{\Delta}t\over 2},t+{{\Delta}t\over 2}] } \over {\Delta}t},
%\end{equation}
%where ${\Delta}t$ can be flexibly determined to fit the dynamism of the intended demand-supply scenario, i.e.,
%to characterize dynamism of the object without considering overly high stochastic volatility.
%We simply call $p$ \emph{power} hereinafter.
As electricity consumption can be investigated on continuous time span of any length,
time domain of the model is described as an arbitrary time interval $\mathcal{T}=[t_1,t_2]$ where $t_1<t_2$ holds permanently.
A function of electric load is a continuous observation of active power defined over time domain: $l(t), t\in\mathcal{T}$.
%characterizes the entire process of an electricity consumption over $\mathcal{T}$.
The entire set of load functions constitute a space
\begin{equation}
\mathcal{L}=\Big\{l(t),t\in\mathcal{T}\Big\},
\end{equation}
namely the space of load functions. Let ${\rm dim}(\cdot)$ denote the dimension of a space.
Next, we translate physical properties of load into algebraic expressions.

\emph{Property 1:}
The addition of parallel connected loads satisfies the following:
\begin{enumerate}
\item  For two independent loads $l_1(t)\in \mathcal{L}$ and $l_2(t) \in \mathcal{L}$, a unique load function is determined by $l_1(t)+l_2(t)=l_3(t)\in\mathcal{L}$.
\item  $l_1(t)+l_2(t)=l_2(t)+l_1(t)$.
\item  $\left(l_1(t)+l_2(t)\right)+l_3(t)=l_1(t)+\left(l_2(t)+l_3(t)\right)$.
\item There exists an element, denoted by $0(t)$ such that $\forall l(t)\in \mathcal{L}, 0(t)+l(t)=l(t)$.
\item $\forall l(t)\in \mathcal{L}$, there exists an element, denoted by $-l(t)$, such that $l(t)+(-l(t))=0$.
\end{enumerate}

\emph{Property 2:}
The scalar multiplication of load functions satisfies the following:
\begin{enumerate}
\item  $\forall a \in \mathbb{R}$ and $\forall l(t)\in \mathcal{L}$, $al(t)\in\mathcal{L}$.
\item  $\forall a_1, a_2 \in \mathbb{R}$ and $\forall l(t)\in \mathcal{L}$, $a_1(a_2l(t))=(a_1a_2)l(t)$.
\item  $(a_1+a_2)l(t)=a_1l(t)+a_2l(t)$.
\item  $a(l_1(t)+l_2(t))=al_1(t)+al_2(t)$.
\item  $1\cdot l(t)=l(t)$.
\end{enumerate}

\emph{Property 3:} $\forall l(t) \in \mathcal{L},$ $\int_{t_1}^{t_2}l^2(t)dt$ exists and is finite. It is in accordance with that load functions are continuous and demander consumes finite amount of energy.

Properties 1 to 3 reveal that:

\emph{Theorem 1:} Load functions constitute a linear space of real-valued functions square integrable on $[t_1, t_2]$~\cite{EO-1989}, i.e.,
\begin{equation}\label{set}
\mathcal{L}\in L_2[t_1,t_2].
\end{equation}

Next, definitions and properties of $L_2$ space~\cite{EO-1989} are presented, which can be applied to the space of load functions $\mathcal{L}$ for further understanding of its physical properties.

\vspace{-1em}
\subsection{Space of Electric Load as $L_2$ Space}

\emph{Definition 1:} For $L_2[t_1,t_2]$, the inner product of $l_1(t)\in L_2[t_1,t_2]$ and $l_2(t) \in L_2[t_1,t_2]$ is defined as
\begin{equation}\label{inner}
\Big(l_1(t),l_2(t)\Big)=\int_{t_1}^{t_2}l_1(t)l_2(t)dt.
\end{equation}
It is an operation to differentiate two arbitrary elements and describe the relationship between the two.
It is the basis to define metrics such as length, distance and orthogonality in a linear space. The inner product on $L_2[t_1,t_2]$ obeys the following:
\begin{enumerate}
\item  $\Big(l_1(t),l_2(t)\Big)=\Big(l_2(t),l_1(t)\Big)$.
\item  $\Big(l_1(t)+l_2(t),l_3(t)\Big)=\Big(l_1(t),l_3(t)\Big)+\Big(l_2(t),l_3(t)\Big)$.
\item  $\forall a\in \mathbb{R}, \Big(al_1(t),l_2(t)\Big)=a\Big(l_1(t),l_2(t)\Big)$.
\item  $\Big(l(t),l(t)\Big)\ge 0$. $\Big(l(t),l(t)\Big)=0 \iff l(t)=0$.
\end{enumerate}

\emph{Definition 2:} The \emph{distance} between two load functions $l_1(t)\in L_2[t_1,t_2]$ and $l_2(t) \in L_2[t_1,t_2]$ is defined as:
\begin{equation}\label{dis}
\rho\Big(l_1(t),l_2(t)\Big)=||l_1(t)-l_2(t)||=\sqrt{\int_{t_1}^{t_2}(l_1(t)-l_2(t))^2dt}.
\end{equation}
By providing a way to quantify the difference between two load functions, it emphasizes the uniqueness of every load compared to one another.

%
%\begin{equation}
%c_k=\int_{t_1}^{t_2}\phi_k(t)p(t)dt.
%\end{equation}

\emph{Definition 3:} The \emph{norm} of a load functions $l(t)$ is defined as:
\begin{equation}
||l(t)||=\sqrt{\Big(l(t),l(t)\Big)}=\sqrt{\int_{t_1}^{t_2}l^2(t)dt}.
\end{equation}
Norm, as can be referred to as a kind of length, provides a way to measure a load function.

%\emph{Definition 3:} The \emph{distance} between two load curves $p_1(t)\in \mathcal{L}$ and $p_2(t) \in \mathcal{L}$ is defined as:
%\begin{equation}
%\rho\Big(p_1(t),p_2(t)\Big)=||p_1(t)-p_2(t)||=\sqrt{\int_{t_1}^{t_2}(p_1(t)-p_2(t))^2dt}.
%\end{equation}
%
%\emph{Property 4:} $L_2[t_1,t_2]$ is a Hilbert space, i.e., an inner product space that is also a complete metric space w.r.t. the distance induced by the inner product.~\footnote{
%For a comprehensive understanding, in Appendix~\ref{AP1}, the operation of inner product is introduced by \emph{Property A.4}; the distance induced by the inner product is introduced by \emph{Definition A.1}; and the properties of a Hilbert space are introduced by \emph{Property A.5}.}

\emph{Property 4:} $L_2[t_1,t_2]$ is a Hilbert space, i.e., an inner product space that is also a complete metric space w.r.t. the distance induced by the inner product.
It satisfies the following:
\begin{enumerate}
\item  $L_2[t_1,t_2]$ is an inner product space based on \emph{Properties 1 to 3} and \emph{Definition 1}.
\item  $L_2[t_1,t_2]$ is a complete metric space w.r.t. $\rho\Big(l_1(t),l_2(t)\Big)$.
\item  $L_2[t_1,t_2]$ is a separable space.
\item  $L_2[t_1,t_2]$ is an infinite dimensional space, i.e., ${\rm dim}(L_2[t_1,t_2])=\infty$.
\end{enumerate}

%\textbf{Remark.} Expression~(\ref{set}) indicates that $\mathcal{L}$ is a subset of $L_2[t_1,t_2]$ as $L_2[t_1,t_2]$ contains elements that can be constructed in math but do not exist in the material world. But all the above definitions and properties can be applied to $\mathcal{L}$.

As the general model, $\mathcal{L}$ is an equivalent description of the set of loads during $\mathcal{T}$ in the real world, where a load function is the unique description of an electricity consumption phenomenon described by energy exchange w.r.t. its existence over $[t_1,t_2]$. The general model has complete capability of distinguishing every load from one another, which is the necessary condition to the goal of modeling, i.e., to perfectly identify every load, clarify the consequence of every load and precisely impact every load. That is the capability of modeling.

$\mathcal{L}$ records complete information w.r.t. electricity consumption activities during $\mathcal{T}$. ${\rm dim}(\mathcal{L})=\infty$ suggests $\mathcal{L}$ and $l(t)$ contains information of infinite amount and diversity, which also reveals the complexity of complete information required to describe loads.
This seems to reveal a negative fact that we may never completely understand a load since we will never be able to capture complete information it delivers. But it is indeed inspiring because:
\begin{enumerate}
\item  Since all information contained in load function is not of equivalent importance w.r.t. reflecting laws of the its nature, it is sufficiency rather than completeness that matters in capturing information and modeling.
\item  To sufficiently understand a load, human can capture as much information as possible by improving observation.
\item  Since the general model explains how a load exists and distinguishes itself from one another, we can use it as a reference and the upper bound to assess the capability of any other load models in capturing information and modeling. The capability of modeling refers to that of distinguishing an object from one another of the same kind.
\end{enumerate}

\emph{Definition 4:} The capability of modeling is defined as the dimension of model space ${\rm dim}(\cdot)$.

\vspace{-0.5em}
\subsection{Analysis of Load Modeling}
In problems such as economic dispatch, unit commitment and power flow~\cite{PG-2013}, load is seldom described as a course of electricity consumption, but an averaged state of the course.
This paradigm of load modeling can be mathematically expressed as:
\begin{equation}
p={1\over t_2-t_1}\int_{t_1}^{t_2}l(t)dt,
\end{equation}
where $p\in \mathbb{R}$ and $[t_1, t_2]$ refers to the intended time interval of which the length usually can be set as 1 minute, 5 minutes, 30 minutes or 1 hour. Let $\mathcal{P}$ denote the space of $p$.
A serial of averaged states are used to describe electricity consumption over a longer time course.
Similarly, in problems of electricity pricing, load is described as a combination of serial averaged states to profile the course of electricity consumption. For example, load profile of a day is described by a vector of hourly energy consumed:
$(E_1,E_2,...,E_{24})$, where $E_k=\int_{T_k}l(t)dt$ and $T_k$ refers to the $k$-th hour of the day, $k=1,2,...,24$. Though the averaged state here is represented by energy instead of power, $E_k$ and $p$ can be regarded as models of the same kind mathematically since they are both defined by the integral of $l(t)$ and the value of time length is a constant once time interval is determined. We take $p$ and $\mathcal{P}$ as the representative to study.

The relationship between $p$ and $l(t)$ can be described by a mapping $\varphi$:
\begin{equation}\label{pp}
p=\varphi\left(l(t)\right):\mathbb{R}^{{\rm dim}(\mathcal{L}_2)}\rightarrow\mathbb{R}^{{\rm dim}(\mathcal{P})}.
\end{equation}
The mapping significantly impacts the modeling capability of $p$ and $\mathcal{P}$ as information w.r.t. elements and relationships in the original space may be lost or distorted during transfer to the image space.
We can assess the capability of $\varphi$ in carrying as many dimensions as possible from the original space to the image space by its rank.
We generalize the definition of rank for linear mapping to general mapping.

\emph{Definition 5:} Rank of mapping $\varphi(a)=b:\mathbb{R}^{{\rm dim}(\mathcal{A})}\rightarrow\mathbb{R}^{{\rm dim}(\mathcal{B})}$ is defined by:
\begin{equation}\label{rank}
{\rm rank}\;\varphi={\rm dim}(\mathcal{B}).
\end{equation}

%Those mappings constitute a space of mappings named the space of cost characteristic, which is denoted by $\mathcal{M}(\mathcal{L}_g,\mathcal{C})$.
%We assess the complexity of $\mathcal{M}(\mathcal{L}_g,\mathcal{C})$ by its dimension. We generalize the definition of dimension for spaces of linear mappings to spaces of general mappings.
%
%\emph{Definition 7:} Given spaces $\mathcal{A}, \mathcal{B}, \mathcal{M}(\mathcal{A},\mathcal{B})$. An element in $\mathcal{M}(\mathcal{A},\mathcal{B}))$ maps an element in $\mathcal{A}$ to an element in $\mathcal{B}$. The dimension of $\mathcal{M}(\mathcal{A},\mathcal{B})$ is defined by:
%\begin{equation}\label{map}
%{\rm dim}(\mathcal{M}(\mathcal{A},\mathcal{B}))={\rm dim}(\mathcal{A}){\rm dim}(\mathcal{B}).
%\end{equation}
%
%Hence, the rank of $\varphi:\mathcal{L}_g\rightarrow\mathcal{C}$ is ${\rm dim}(\mathcal{C})$. The dimension of $\mathcal{M}(\mathcal{L}_g,\mathcal{C}))$ is ${\rm dim}(\mathcal{L}_g){\rm dim}(\mathcal{C})$.
To analyze the space $\mathcal{P}$, we present the definition of isomorphic space and isomorphic mapping as a preliminary.

\emph{Definition 6:} Two linear spaces $\mathcal{A}$ and $\mathcal{B}$ are isomorphic, if a bijective mapping $f:\mathbb{R}^{{\rm dim}(\mathcal{A})} \to \mathbb{R}^{{\rm dim}(\mathcal{B})}$ satisfies:
$\forall \alpha,\beta \in \mathbb{R}, l_1(t),l_2(t)\in \mathcal{A}$,
\begin{equation}\label{iosom}
f(\alpha l_1+\beta l_2)=\alpha f(l_1)+\beta f(l_2).
\end{equation}
Such a mapping is called an isomorphic mapping.

According to equation~(\ref{pp}), $p$ is a real number, and the space $\mathcal{P}$ is therefore a field of real numbers. We have
\begin{equation}
 {\rm rank}\;\varphi={\rm dim}(\mathcal{P})=1.
\end{equation}
Obviously, a continuous load function has been compressed to a real number, and the dimension of the original space has been dramatically reduced to 1 from $\infty$.
It is like we use a piece of continuously exposed film of starry sky to conclude the course of changes in the sky that night.
The modeling capability of $\mathcal{P}$ is 1, which indicates that $\mathcal{P}$ is far from being an isomorphic space of $\mathcal{L}$.
Physical explanation of the mathematical form is that
a considerably amount of information w.r.t. elements and relationships in the set of real electric loads is lost and the remained information carried by the sole dimension may be distorted. We use a case to further clarify the capability of modeling of this paradigm.

\begin{figure}[!t]
\centering
%\vspace{-0.5em}
\includegraphics[width=\linewidth]{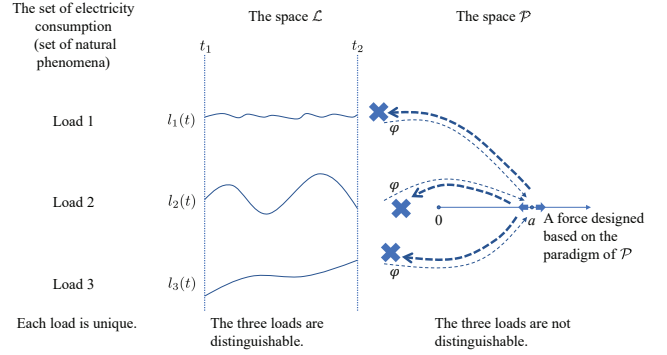}
%\vspace{-1.5em}
\caption{A case of loads with equal averaged power.} \label{case1}
%\vspace{-1em}
\end{figure}

{\bf\emph{Case 1:}} Consider three loads in Fig.~\ref{case1}. As three independent electricity consumption courses, they are respectively and uniquely characterized by three load functions $l_1(t), l_2(t), l_3(t)$. They satisfy $\int_{t_1}^{t_2}l_1(t)dt=\int_{t_1}^{t_2}l_2(t)dt=\int_{t_1}^{t_2}l_3(t)dt$.
In the space $\mathcal{L}$, it is easy to distinguish each of them by checking the form of load functions, calculating inner product or distance between arbitrary two according to formula~(\ref{inner}) or formula~(\ref{dis}). The mapping $\varphi$ has respectively mapped $l_1(t), l_2(t), l_3(t)$ into the same point $a={1\over t_2-t_1}\int_{t_1}^{t_2}l_1(t)dt$ in the space $\mathcal{P}$. This means in the space $\mathcal{P}$, the three loads are not distinguishable. If a force to reshape load is designed based on the modeling paradigm of $\mathcal{P}$, i.e., defined in the space of $\mathcal{P}$, it can only impact the value of $a$ in the field of real numbers. But the force can never be conducted from $\mathcal{P}$ to $\mathcal{L}$ to reshape the form of load since $\varphi$ is not a bijective mapping. Actually, $\{l(t)|{1\over t_2-t_1}\int_{t_1}^{t_2}l(t)dt=a\}$ refers to a hyperplane of $({\rm dim}(\mathcal{L})-1)$ dimensions that belongs to $\mathcal{L}$.
The capability of the force in distinguishing load is limited to distinguishing hyperplanes in $\mathcal{L}$, while all loads in the hyperplane impacted are only forced to adjust the integral of load function instead of the form. The goal of reshaping the form of load is far from being achieved.

%\vspace{-0.5em}
\section{Electricity Pricing: General Model and Analysis of Classic Modeling}\label{A2}

\subsection{The General Model}
Let $\rho\in\mathbb{R}$ denote the general payment of an arbitrary load function $l(t)\in\mathcal{L}$ determined by a mapping $\lambda:\mathbb{R}^{{\rm dim}(\mathcal{L})}\rightarrow\mathbb{R}$, also known as the pricing. A scalar field
\begin{equation}\label{ori}
\rho=\lambda\left(l(t)\right)
\end{equation}
describes the relationship between payment and load function.
It is a hyperplane of ${\rm dim}(\mathcal{L})$ dimensions that belongs to a space of $({\rm dim}(\mathcal{L})+1)$ dimensions.
Gradient $\nabla \rho$ indicates the direction that payment increases at the fastest rate. Therefore, $(-\nabla \rho)$ represents the direction of incentive or force of payment to reshape the form of a load function with maximum payment reduced.

Payment of electricity generally have to address two concerns: to reflect the cost of electricity supply
%\footnote{As a primary theoretical research to investigate the relationship among electricity payment, generation cost and electricity consumption, it does not consider the cost of transmission.}
and to impact electricity consumption, i.e. the shape of load.
Therefore, capability of modeling of the general pricing model must not be lower than either that of the electricity supply cost model or that of the general load model.

{\bf\emph{Proposition 1:}} The general pricing model has a complete capability in shaping the form of a load.
\begin{proof}
$\rho=\lambda\left(l(t)\right)$ is defined over the space $\mathcal{L}$. The force $(-\nabla \rho)$, as is of ${\rm dim}(\mathcal{L})$ dimensions, is composed of partial derivatives of all dimensions of $\mathcal{L}$, which refers to a force of full capability in impacting every dimension of $l(t)$.
\end{proof}

Electricity supply cost, as a consequence of electricity consumption, can be described by a function of load.
Let $\gamma\in\mathbb{R}$ denote the general supply cost of an arbitrary load function $l(t)\in\mathcal{L}$ determined by a mapping $\varphi:\mathbb{R}^{{\rm dim}(\mathcal{L})}\rightarrow\mathbb{R}$. A scalar field
\begin{equation}\label{cost}
\gamma=\varphi\left(l(t)\right)
\end{equation}
describes the relationship between supply cost and load function.
It is also a hyperplane of ${\rm dim}(\mathcal{L})$ dimensions that belongs to a space of $({\rm dim}(\mathcal{L})+1)$ dimensions.
Gradient $\nabla \gamma$ indicates the direction that supply cost increases at the fastest rate, i.e., cost characteristic of the system as its inherent property.
Therefore, $(-\nabla \gamma)$ refers to the direction that supply cost decreases at the fastest rate.

{\bf\emph{Proposition 2:}} The best pricing to minimize the supply cost is to let $\nabla \rho= a\nabla \gamma$, where $a>0$ is a constant.

\begin{proof}
$-\nabla \rho= -a\nabla \gamma$ refers to that the direction of incentive or force of payment to reshape the form of a load function with maximum payment reduced becomes exactly the one that supply cost decreases at the fastest rate.
\end{proof}

{\emph{Proposition 1}} and {\emph{Proposition 2}} clarify law of interaction between payment and load at the demand side and the necessity of pricing to sufficiently convey information on how load results in cost.

\subsection{Analysis of Classic Electricity Pricing}
The classic pricing is integral-based, as the payment can be calculated as integral of power times unit price of energy or power capacity\footnote{A kind of averaged power that can be  written as integral of power divided by time interval.} times unit price of power capacity.
Here we choose the former to investigate as similar judgement can be made on the latter. The investigation is also conducted over time interval $\mathcal{T}$.
Let $\rho_{\rm{cls}} \in \mathbb{R}$ denote the classic payment, $\mathcal{R}_{\rm{cls}}$ the set of classic payments and $\lambda_{\rm{cls}}$ the classic pricing, respectively.
We have
\begin{equation}\label{cls}
\rho_{\rm{cls}}=\lambda_{\rm{cls}}(\int_{t_1}^{t_2} l(t) {\rm d}t)=\rho_{\rm{unit}}\int_{t_1}^{t_2} l(t) {\rm d}t,
\end{equation}
where $\rho_{\rm{unit}}$ is a constant unit price predetermined by
\begin{equation}
\rho_{\rm{unit}} = {F \over \int_{t_1}^{t_2} g(t) {\rm d}t}.
\end{equation}
$\int_{t_1}^{t_2} g(t) {\rm d}t$ refers to the gross energy delivered and $F$ the gross supply cost.
The mapping $\lambda_{\rm{cls}}$ actually maps the integral of a load function, instead of $l(t)$ itself, to a number.
That is, the classic pricing is defined over the compressed space $\mathcal{P}$ instead of $\mathcal{L}$.
We have ${\rm{dim}}(\mathcal{P})=1$, ${\rm{dim}}(\mathcal{R}_{\rm{cls}})=1$.

{\bf\emph{Proposition 3:}} Classic integral-based pricing model is ineffective in reshaping load over an arbitrary interval $\mathcal{T}$.
\begin{proof}
The pricing $\lambda_{\rm{cls}}: \mathbb{R}\rightarrow \mathbb{R}$ is defined over $\mathcal{P}$ instead of $\mathcal{L}$.
Dimension of the hyperplane defined by $\rho_{\rm{cls}}=\lambda_{\rm{cls}}(\int_{t_1}^{t_2} l(t){\rm d}t)$ is $1$ compared to
the dimension of the original hyperplane defined by $\rho=\lambda\left(l(t)\right)$ which is $\infty$, indicating that
the original hyperplane has be degraded to a line.
We further investigate how the incentive impact elements of the original space.
An element $a \in \mathcal{P}$ refers to a hyperplane of $\mathcal{L}$: $\{l(t)|{1\over t_2-t_1}\int_{t_1}^{t_2}l(t)dt=a\}$.
The impacting force of payment can be described by
\begin{equation}
 -{ {\rm d} \rho_{\rm{cls}} \over {\rm d} a} = -\rho_{\rm{unit}},
\end{equation}
where $\rho_{\rm{unit}}$ refers to the absolute strength and the symbol $(-)$ refers to an incentive to decrease $a$.
Its impact is on a hyperplane instead of an element of $\mathcal{L}$ that fails to distinguish loads in $\mathcal{L}$.
Consequently, the impact is on the integral of loads in the same hyperplane instead of the form of loads, which fails to reshape the form of a load.
The result is uncertain reshaping of every load and the aggregated load over $\mathcal{T}$.
\end{proof}
It is more common that eq.~(\ref{cls}) is written as $\rho_{\rm{cls}}=\lambda_{\rm{cls}}(p)$, which is independent of time. But the rigorous expression actually should be
\begin{equation}
\rho_{\rm{cls}}=\lambda_{\rm{cls}}(p), \forall t,
\end{equation}
where $\forall t$ reveal the premise of application is that we assume the relationship holds at any moment.

A common solution to improve the capability of pricing in impacting load is to shorten the length of pricing cycles, which is known as the spot pricing~\cite{OP-1982,EW-2003}.
The spot price for buying and selling electric energy is determined by the supply and demand conditions at that instant. The term \emph{instant} refers to a cycle $\Delta t$ that is much shorter than the time scale of the studied phenomena $\mathcal{T}$, e.g., 5 minutes compared to 1 day, in which physical properties of the studied objects is assumed to be fixed for modeling. In another word, phenomena and relationships in a cycle are described by static model. The \emph{term} static refers to the classic paradigm of modeling that describe the relationship between variables independent of time and time course, i.e, $f$ in $a=f(b)$ that permanently describes the relationship between $a$ and $b$ despite that $b$ varies with time. Actually, the rigorous expression is $a=f(b), \forall t$, in which we usually ignore $\forall t$. But the condition $\forall t$ strongly restricts the rank of $f$. $a=f(b), \forall t$ is totally different from $a=f(b(t)), t\in[t_1,t_2]$ or $a(t)=f(b(t)), t\in[t_1,t_2]$.
Examples are the locational marginal price~(LMP), real-time pricing and hourly spot pricing.
For $\Delta t$ and $\mathcal{T}$ we have $N\Delta t=t_2-t_1, N\in \mathbb{N}_+$.
A series of spot unit prices $((\rho_{\rm{unit}})_1,(\rho_{\rm{unit}})_2,...,(\rho_{\rm{unit}})_N)$ is defined for the $N$ cycles $\{[t_1,t_1+\Delta t),[t_1+\Delta t, t_1+2\Delta t),...,[t_2-\Delta t,t_2)\}$.
The payment over $\mathcal{T}$ is expressed as:
\begin{equation}\label{totalp}
\rho_{\rm{cls}}=\sum_{k=1}^N (\rho_{\rm{unit}})_k\int_{t_1+(k-1)\Delta t}^{t_1+k\Delta{t}} l(t) {\rm d}t.
\end{equation}
Consider the $N$ integrals $\{\int_{t_1+(k-1)\Delta t}^{t_1+k\Delta{t}} l(t) {\rm d}t; k=1,2,...,N \}$ as $N$ independent variables, we can realize an $N$-dimensional hyperplane defined by eq.~(\ref{totalp}) compared to the $1$-dimensional hyperplane defined by eq.~(\ref{cls}).
It appears that the capability of pricing in distinguishing and impacting load is improved by the division of a long cycle into $N$ shorter cycles. But capability of the pricing in reflecting the real supply cost has not been improved. There are three major reasons.
\begin{enumerate}
\item Pricing that is defined over the hyperplane $\{l(t)|\int_{t_1+(k-1)\Delta t}^{t_1+k\Delta{t}} l(t) {\rm d}t=a_k,a_k\in\mathbb{R},k=1,2,...,N \}$ cannot carry the information on relationship between supply cost and load that have to be defined over a time interval of $\mathcal{T}$ of which the length belongs to $(\Delta t, t_2-t_1]$. That is, by dividing a long cycle into multiple shorter ones, we break the continuity of the operation in time domain and cannot retain complete information on the physical properties associated with undivided course of electricity consumption of which the time length is an arbitrary positive number between 0 and $(t_2-t_1)$.
\item Since $\Delta t$ refers to a course of time, $l(t), t\in [0,\Delta t]$ is also a continuous load function that belongs to a space of infinite dimensions. Problem stated in \emph{Proposition 3} still exists in every $\Delta t$.
\item Despite that $N$ cycles creates $N$ independent dimensions and $N$ is varying, all cycles apply a unified, predetermined, static supply cost model that does not utilize any extra information on relationship between supply cost and load which is associated with the time length of $\Delta t$. The magnitude of error between the classic model and the reality varies over the specification of time length $\Delta t$, e.g., 5 minutes, 1 minute or 30 seconds, which is uncontrollable. With the decrease in $\Delta t$, it becomes more and more difficult for the static model to maintain a consistent validity in every $\Delta t$.
\end{enumerate}
We present a case to clarify the problem.

{\bf\emph{Case 2:}} LMP refers to the marginal cost of supplying the next increment of electric energy at a specific bus considering the generation marginal cost and the physical aspects of the transmission system~\cite{MOE-2002,OP-1982}. The LMP is defined as
\begin{equation}
{{\rm{LMP}}_k}= {{\rm{LMP}}_k^{\rm{ref}}}+{{\rm{LMP}}_k^{\rm{loss}}}+{{\rm{LMP}}_k^{\rm{cong}}},
\end{equation}
where ${{\rm{LMP}}_k^{\rm{ref}}}$ is marginal generation cost at the reference bus, ${{\rm{LMP}}_k^{\rm{loss}}}$ is marginal transmission loss and ${{\rm{LMP}}_k^{\rm{cong}}}$ marginal congestion price. We take ${{\rm{LMP}}_k^{\rm{ref}}}$ as an example. The determination of ${{\rm{LMP}}_k^{\rm{ref}}}$ depends on cost characteristics of the related generation units, which are usually modeled by $F$-$P$ curves. The $F$-$P$ curve is a fitted curve based on observations on quasi-static operation of generator~\cite{PSE-2018,PG-2013}. With the quasi-static operation, the only 2 dimensions and the operation of fitting, the $F$-$P$ curve contains extremely limited information on cost characteristic of the generation unit while the remained information is possibly distorted. Despite that the specification of time length $\Delta t$ varies and dynamism and uncertainty of system operation in each $\Delta t$ is different, ${{\rm{LMP}}_k^{\rm{ref}}}$ is calculated based on the same generation cost characteristic models in all cases. It is with no doubt the ${{\rm{LMP}}_k^{\rm{ref}}}$ considerably deviates from the real generation costs.

It is not sufficient to solely improve the capability of pricing in distinguishing and impacting load. What the classic modeling paradigm has not achieved is to create and exploit available dimensions of model to capture as much as possible the key information that describes the physical properties of studied objects, e.g., the relationship between supply cost and load function. To do this, one has to be fully conscious that interactions among objects exist in the form of undivided time course. Models of the corresponding properties and relationships have to be defined on a time course instead of being independent of time.

%\vspace{-0.5em}
\section{Specialized Models and Examples}\label{PRPS}
Eq.~(\ref{cost}) has fully abstracted the relationship between supply cost and load function, and eq.~(\ref{ori}) the relationship between payment and load function.
But physical properties cannot be further perceived from those abstract forms, let alone applications. Computable model as specific form of abstract model is in need. Yet it is hard to establish computable mapping from the original space $\mathcal{L}$ to the set of supply cost or that of payment.

In this section, we introduce isomorphic vector space of the original space $\mathcal{L}$ that is constituted by orthonormal basis, where physical properties of load can be described in analytical forms and computable mappings can be defined. On the isomorphic space, we present general form of the proposed computable supply cost model and pricing model.
Finally, an applicable model based on the Fourier Series is presented.
%\emph{Property A.6:} Power curves $p(t)\in \mathcal{L}$ satisfy the Dirichlet conditions:
%\begin{enumerate}
%\item $p(t)$ has a finite number of finite discontinuities.
%\item $p(t)$ has a finite number of maxima and minima.
%\item  $\int_{t_1}^{t_2}|p(t)|dt$ is finite.
%\end{enumerate}
\vspace{-0.5em}
\subsection{Isomorphic Space of $\mathcal{L}$ }
Since $\mathcal{L}$ is a Hilbert space, elements in $\mathcal{L}$ can be decomposed into a series of coefficient-weighted orthonormal basis functions.

\emph{Definition 7:} An orthonormal basis is an infinite set of functions that satisfy:
for $\phi_i(t)$ and $\phi_j(t)$ of the set,
\begin{equation}
\Big(\phi_i,\phi_j\Big) = \left\{ \begin{array}{ll}
1   & \text{if } i = j,\\
0  & \text{else}.
\end{array} \right.
\end{equation}

The space $\mathcal{L}$ contains many orthonormal bases such that:

\emph{Property 5:} A load curve $l(t)$ in $\mathcal{L}$ can be expressed as a linear combination of the functions of an orthonormal basis $\{ \phi_k \}$ as:
\begin{equation}\label{series}
l(t)=\sum_k \mu_k\phi_k,
\end{equation}
where $\mu_k=\Big(l(t),\phi_k\Big)$.

\emph{Property 6:} In equation~(\ref{series}), $l(t)$ and $\{\mu_k\}$ always satisfy the Parseval's theorem: $||l(t)||^2=\sum_k \mu_k^2$.

We define a vector $\bm l=(\mu_1,\mu_2,...)$, between which and $l(t)$ an isomorphic mapping always exists such that each $l(t)$ is compared to a point $\bm l$ in a space of multiple dimensions. All the $\bm l$ constitute a vector space $\hat{\mathcal{L}}$ that is isomorphic to $\mathcal{L}$.
We can compare a load in the real world to a point in the orthonormal vector space $\hat{\mathcal{L}}$.
$\bm l \in \hat{\mathcal{L}}$, just as $l(t)\in\mathcal{L}$, distinguishes a load from one another.

In the vector space $\hat{\mathcal{L}}$, the capability of modeling can be interpreted as to identify every point in accordance with reality. For electricity pricing, to reshape the form of load can be interpreted as to ensure every point is identified and imposed with precise economic force towards right direction.
Ineffectiveness of the classic modeling is interpreted as incapability in distinguishing a point from one another in $\hat{\mathcal{L}}$.
The integral-based pricing cannot identify every point, such that the economic force impacts all points that belongs a same hyperplane while the force does not indicate a specific direction for the movement of points. The movement of point is not strictly regulated, which explains the uncertain reshaping of every load.

\vspace{-1em}
\subsection{The Specialized Algebraic Model}
\subsubsection{Space of Dynamism} Since load curve reveals the rhythm of human activities, we are interested in the isomorphic orthonormal vector space of $\mathcal{L}$ whose basis functions have rhythmic dynamism, i.e., being periodic. That is, the space is constituted by orthonormal basis with definable and distinguishable periodicity, and is simply named the space of dynamism.
Let $T_0=t_2-t_1$ denote length of the period. Let $f_0= {1\over T_0}$, which is referred to as the fundamental frequency over interval $\mathcal{T}$. We intend to use an orthonormal basis $\{ \phi_k|k=0,1,2,... \}$ that satisfy the following:

\emph{Property 7:} $\phi_0$ is named the zero-frequency basis function, which satisfies
\begin{equation}
\phi_0(t)\equiv {1\over T_0}.
\end{equation}
Obviously, $\phi_0(t)$ has no dynamism and
\begin{equation}
\mu_0={1\over T_0}\int_{t_1}^{t_2} l(t)dt.
\end{equation}
That is, $\phi_0$ can independently characterize energy consumption of $l(t)$ with respect to accumulation of energy consumed.
$\phi_k, k > 0$ satisfy
\begin{equation}
\phi_k(t+{T_0\over k})=\phi_k(t).
\end{equation}
We say that $\phi_k$ has dynamism as it fluctuates with frequency $kf_0$. We further require that
\begin{equation}
\int_{t_1}^{t_2} \phi_k(t)dt=0.
\end{equation}
That is, the basis functions with dynamism accumulate zero energy in time interval $\mathcal{T}$. Thus, dynamism is independent of energy accumulation that is solely characterized by $\phi_0$.

For $\phi_i$ satisfying $\phi_i(t+{T_0\over i})=\phi_i(t)$ and $\phi_j$ satisfying $\phi_j(t+{T_0\over j})=\phi_j(t)$, if
\begin{equation}
i>j,
\end{equation}
we say that the dynamism of $\phi_i$ with regard to frequency is higher than that of $\phi_j$.

A load curve can thus be decomposed by $\{ \phi_k|k=0,1,2,... \}$ into two major parts: the accumulative part measured by $T_0 c_0 \phi_0$, the dynamic part described by the other basis functions $\{ \phi_k|k=1,2,... \}$. Moreover, the dynamism of a load curve has to be analyzed from more than just the aspect of frequency.
The other aspect is strength, i.e., coefficients of the basis $\{\mu_k,k\in \mathbb{Z}\}$ indicate the polarity and amplitude of frequency component $\phi_k$.
${\bm l}=(\mu_0,\mu_1,...)$ denotes a vector representing a power function. The space is denoted by $\hat{\mathcal{L}}$.

%A one-to-one mapping $\mathcal{F}$ can be defined between $\bf c$ and $p(t)$ as
%\begin{equation}
%{\bf c}=\mathcal{F}\left(p(t)\right): \mathcal{L} \to \mathcal{H}.
%\end{equation}
%$\mathcal{H}$ is a dual equivalent of $\mathcal{L}$.
%As the vector can fully characterize the dynamism of load, it is named the vector of load dynamism. $\mathcal{H}$ is therefore named the space of load dynamism.
\subsubsection{The Supply Cost and Cost Characteristic}
In the isomorphic space of dynamism $\hat{\mathcal{L}}$, eq.~(\ref{cost}) can be rewritten as
\begin{equation}\label{cost2}
\gamma=\varphi({\bm l}).
\end{equation}
$\nabla \gamma$ can be derived as:
\begin{equation}
\nabla \gamma =\left({\partial \gamma\over \partial \mu_0},{\partial \gamma\over \partial \mu_1},...\right).
\end{equation}
$\nabla \gamma$ has sufficient dimensions, i.e., equivalent dimensions that $\bm l$ has, to completely describe the cost characteristic of the system.
${\partial \gamma\over \partial \mu_k}=\iota_k \in\mathbb{R}$ is a scalar that refers to an aspect~(a dimension) of the cost characteristic of the system, and can be calculated via sufficient observations. We have $\nabla \gamma=(\iota_0,\iota_1,...)$.
We can let
\begin{equation}\label{ori3}
\gamma=\sum_{k\ge 0} \iota_k\mu_k,
\end{equation}
which has a simple form with complete capability of modeling, i.e., ${\rm dim}(\nabla \gamma)={\rm dim}(\mathcal{L})$.

\subsubsection{The Payment and Pricing}
In $\hat{\mathcal{L}}$, eq.~(\ref{ori}) can be rewritten as
\begin{equation}\label{ori2}
\rho=\lambda({\bm l}).
\end{equation}
$\nabla \rho$ can be derived as:
\begin{equation}
\nabla \rho =\left({\partial \rho\over \partial \mu_0},{\partial \rho\over \partial \mu_1},...\right).
\end{equation}
$\nabla \rho$ has sufficient dimensions, i.e., equivalent dimensions that $\bm l$ has.
${\partial \rho\over \partial \mu_k}=\lambda_k \in\mathbb{R}$ is a scalar that refers to a rate~(a dimension) of the pricing rates. We have $\nabla \rho=(\lambda_0,\lambda_1,...)$.
We can let
\begin{equation}\label{ori3}
\rho=\sum_{k\ge 0} \lambda_k\mu_k,
\end{equation}
which has a simple form with complete capability of modeling, i.e., ${\rm dim}(\nabla \rho)={\rm dim}(\mathcal{L})$.
According to \emph{Proposition 2},
we determine the payment by
\begin{equation}\label{total}
\rho= T_0 \sum_{k\ge 0} \lambda_k\mu_k= T_0 \sum_{k\ge 0} \iota_k\mu_k.
\end{equation}

In~(\ref{total}), $\lambda_0$ and $T_0\lambda_0\mu_0$ can be regarded as the classic unit price and the total payment for the amount of energy consumed, respectively. $\lambda_k,k>0$ is the unit price for adding dynamism to the system described by $\mu_k\phi_k,k>0$.
This suggests that the proposed pricing can be regarded as a generalization of the classic pricing.

\vspace{-1em}
\subsection{Computational Model Based on the Fourier Series}\label{PMFS}
The Fourier series~\cite{EO-1989} is a specific orthonormal basis that is easy to compute and has been widely applied to analyze electrical phenomena. Moreover, it can distinguish dynamism of a load curve from the energy accumulated as the Fourier series satisfy \emph{Property 7}.
We introduce the Fourier Series as a specific case of $\{ \phi_k|k=0,1,... \}$.
Based on them we further present a pricing method.

\subsubsection{The Fourier Series}
%\subsection{The Fourier Series and Fourier Transform}
In $L_2[t_1,t_2]$, a set of functions
\begin{small}
\begin{equation}\nonumber
{1 \over T_0}, \cos {2\pi\over T_0} t, \sin {2\pi \over T_0} t,\cos {4\pi \over T_0} t, \sin {4\pi \over T_0} t,...\cos 2\pi{n \over T_0} t, \sin 2\pi {n\over T_0} t,...
\end{equation}
\end{small}
constitute an orthonormal basis. Each of the basis functions is also called a frequency component. ${1 \over T_0}$ is the zero-frequency component. $\cos {2\pi\over T_0} t$ and $\sin {2\pi \over T_0} t$ are the fundamental frequency components. The other components each has an integer multiples of the fundamental frequency.

A load curve $l(t)\in \mathcal{L}$, as it naturally satisfies the Dirichlet conditions~\footnote{The Dirichlet conditions is the sufficient conditions for a real-valued function on an intended interval to be equal to the sum of its Fourier series.}, can be written as:
\begin{equation}\label{FS}
l(t)={a_0\over 2}+\sum_{n=1}^\infty\Big[a_n\cos(2\pi{n f_0}t)+b_n\sin(2\pi{n f_0}t) \Big],
\end{equation}
where $n\in \mathbb{Z}_+$, coefficients
\begin{equation}\label{FS1}
a_n={2\over T_0}\int_{t_1}^{t_2} l(t)\cos(2\pi{n f_0}t)dt
\end{equation}
and
\begin{equation}\label{FS2}
b_n={2\over T_0}\int_{t_1}^{t_2} l(t)\sin(2\pi{n f_0} t)dt
\end{equation}
respectively indicate the polarities and amplitudes of the corresponding sine and cosine components.
The right side of equation~(\ref{FS}) is called the Fourier series. Obviously,
\begin{equation}\label{FS3}
a_0 ={2\over T_0}\int_{t_1}^{t_2} l(t)dt,%~\footnote{Generally, we use $2\over T_0$ instead of $1\over T_0$ as the zero-frequency basis function in a Fourier series.}
\end{equation}
and ${T_0\over 2}a_0$ equals the energy consumed. Hence, the Fourier coefficients $\{(a_n,b_n)|n\in\mathbb{Z}\}$ fully characterize a load curve w.r.t. energy consumed and the dynamism.

%If we generalize $p(t)$ from $[t_1,t_2]$ to $[-\infty,+\infty]$ and further require that $\int_{-\infty}^{+\infty}|p(t)|dt$ is finite, we can apply the Fourier transform to express it as:
%\begin{equation}\label{FT}
%p(t)=\int_0^\infty\Big[a(f)\cos(2\pi ft)+b(f)\sin(2\pi ft) \Big]df,
%\end{equation}
%where $f= {1\over T}$ is referred to as frequency, and the coefficient functions $a(f)$ and $b(f)$, generalized from the Fourier coefficients, are respectively given by:
%\begin{equation}
%a(f)={2}\int_{-\infty}^{+\infty} p(t)\cos(2\pi ft)dt
%\end{equation}
%and
%\begin{equation}
%b(f)={2}\int_{-\infty}^{+\infty}  p(t)\sin(2\pi{f} t)dt.
%\end{equation}

\subsubsection{Price-Frequency Coefficients}
For the zero-frequency component, let $\alpha_0 >0$ denote the price per unit of energy consumed. For the cosine components, let $\alpha_n, n>0$ represent the unit price for $a_n\cos(2\pi{n f_0}t)$. For the sine components, let $\beta_n, n>0$ represent the unit price for $b_n\sin(2\pi{n f_0}t)$.
%Since $\alpha_0$ is defined for energy while $\alpha_n$ and $\beta_n$ are both defined for the Fourier coefficients, the dimension of $\alpha_0$ differs from the other one that is used by both $\alpha_n$ and $\beta_n$.
The payment of a load is defined by
\begin{equation}\label{pay2}
\rho= \alpha_0 {T_0\over 2}a_0 + T_0\sum_{n=1}^{\infty} (\alpha_n a_n+\beta_n b_n).
\end{equation}
Hence, we have
\begin{equation}
\nabla \rho=T_0({\alpha_0\over 2},\alpha_1,\beta_1,\alpha_2,\beta_2,...),
\end{equation}
where $\alpha_n$ and $\beta_n$ actually denote the same objects as $\lambda_k$ does. To determine these coefficients we need to evaluate $\iota_k$.
Since we can observe $\gamma$ and $l(t)$ as many times as possible, by finite number of observations we can formulate linear equations to evaluate $\iota_k$
according to~(\ref{ori3}), (\ref{FS}), (\ref{FS1}), (\ref{FS2}), and (\ref{FS3}).

\begin{figure}[!t]
\centering
%\vspace{-0.3em}
\includegraphics[width=\linewidth]{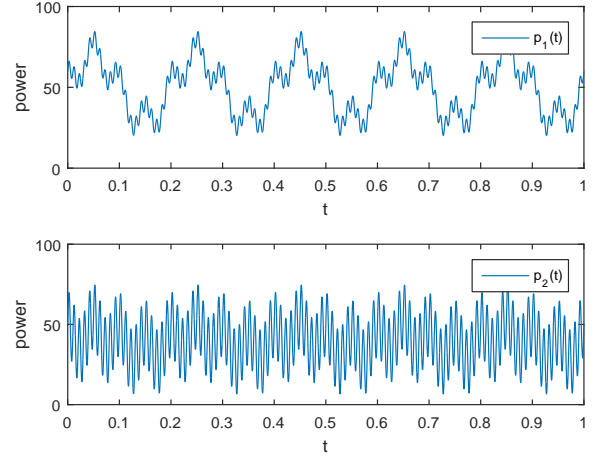}
%\vspace{-1em}
\caption{Load prifles:$l_1(t)=50+20\sin(10\pi t)+10\cos(40\pi t)+5\sin(200\pi t)$, $
l_2(t)=40+5\sin(10\pi t)+10\cos(40\pi t)+20\sin(200\pi t)$.} \label{sload}
%\vspace{-2em}
\end{figure}

\subsection{Examples}\label{EXMP}
Examples based on continuous functions are presented to illustrate principles of the proposed pricing emphasizing conciseness and generality, where practical issues such as units of measurement and physical bases of price-frequency functions are not considered. All load curves are defined over $[0,1]$, which means $T_0=1$.
For comparison, we assess the payments of two load curves that are respectively defined by
\begin{equation}
\begin{aligned}
l_1(t)=50+20\sin(10\pi t)+10\cos(40\pi t)+5\sin(200\pi t),\\
l_2(t)=40+5\sin(10\pi t)+10\cos(40\pi t)+20\sin(200\pi t),
\end{aligned}
\end{equation}
which are drawn in Fig.~\ref{sload}.
For $l_1(t)$, $(a_0)_1=100$, $(b_5)_1=20$, $(a_{20})_1=10$ and $(b_{100})_1=5$.
For $l_2(t)$, $(a_0)_2=80$, $(b_5)_2=5$, $(a_{20})_2=10$ and $(b_{100})_2=20$.
Obviously, Load 2 consumes less energy than Load 1 but is of higher dynamism.

\begin{figure}[!t]
\centering
%\vspace{-0.3em}
\includegraphics[width=\linewidth]{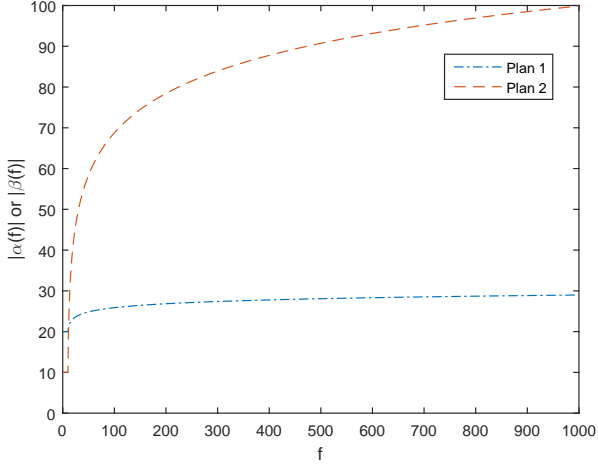}
%\vspace{-1em}
\caption{Price-frequency coefficient functions.} \label{ss}
%\vspace{-2em}
\end{figure}

We also assume two billing plans, i.e. two sets of price-frequency coefficient functions that reflect the cost characteristics of two different supply conditions, respectively.
As shown in Fig.~\ref{ss}, we define $\{[|\alpha(f)|]_1,[|\beta(f)|]_1\}$ for Plan 1 as:
\begin{equation}
\begin{aligned}
\;\;[|\alpha(f)|]_1&=[|\beta(f)|]_1\\
&=\left\{ \begin{array}{ll}
20   & \text{if } 0\le f < 10,\\
3\log_{10}(f-9)+20  & \text{else}.
\end{array} \right.
\end{aligned}
\end{equation}
and
$\{[|\alpha(f)|]_2,[|\beta(f)|]_2\}$ for Plan 2 as:
\begin{equation}
\begin{aligned}
\;\;[|\alpha(f)|]_2&=[|\beta(f)|]_2\\
&=\left\{ \begin{array}{ll}
10   & \text{if } 0\le f < 10,\\
30\log_{10}(f-9)+10  & \text{else}.
\end{array} \right.
\end{aligned}
\end{equation}
They are presented in the form of absolute value.
To determine the final value of a price-frequency coefficient, we assume that the polarity of the price-frequency coefficient should be same with that of the corresponding Fourier coefficient at the supply side, e.g. $\alpha(f)a(f)>0$ and $\beta(f)b(f)>0$.\footnote{The assumption is made for presentation of the basic idea at a theoretical
level. We note that $\alpha(f)a(f)>0$ or $\beta(f)b(f)>0$ may not always hold in practice. But so far no experiment has been conducted to study the relationship between cost-frequency coefficients and the corresponding Fourier coefficients w.r.t. the dynamic characteristics of a generator.}
In a one-source-one-subscriber case, the generation curve equals the load curve. Hence, for Plan 1, $[\alpha(0)]_1=20$, $[\beta(5)]_1=20$, $[\alpha(20)]_1=23.9031$ and $[\beta(100)]_1=26$.
For Plan 2, $[\alpha(0)]_2=10$, $[\beta(5)]_2=10$, $[\alpha(20)]_2=49.0309$ and $[\beta(100)]_1=70$.
Compared to Plan 1, Plan 2 charges for the energy consumed at a lower rate but charges for the dynamism at a much higher rate.
That is, the supply condition of Plan 1 can cope with load of relatively high dynamism at a relatively low supply cost while that of Plan 2 obtains a relatively lower supply cost when the load is of relatively low dynamism.

Bills of the four combinations are presented in TABLE~\ref{paytbl}. Since Load 1 consumes more electricity than Load 2, Load 1 is charged more than Load 2 for the non-dynamic part by either plan. Since Load 2 poses heavier burdens on the system than Load 1 w.r.t. dynamism of load, Load 2 is charged more than Load 1 for the dynamic part by either plan.
It suggests that the supply of Plan 1 is more cost-efficient in serving load of higher dynamism, e.g., Load 2, while that of Plan 2 is more cost-effective in serving load of lower dynamism, e.g., Load~1.

\begin{table}
\caption{Bills of the four combinations}
\label{paytbl}
\centering
\begin{tabular}{|l|c|c|c|}
\hline
\diagbox{Combo.}{Value}{Item} & Non-dyn. & Dyn.& Total \\
\hline
Load1Plan1 & 1000 & 769.031  & 1769.031\\
\hline
Load2Plan1 &  800 & 859.031 & 1659.031 \\
\hline
Load1Plan2 &  500 & 1040.309&  1540.309\\
\hline
Load2Plan2 &  400 & 1940.309 & 2340.309\\
\hline
\end{tabular}
\end{table}

\section{Conclusion}\label{clnc}
By abstractly describing features and relationships of demand-supply activity of power system on time course of any length, we have defined the general models of load, supply cost and pricing that have complete capability in modeling. By describing the classic paradigm of modeling load, supply cost and pricing in a general algebraic sense and analyzing it under the framework of general models with complete modeling capability, we have proved the ineffectiveness of the classic paradigm of modeling and presented sufficient explanation on the causes.
The algebraic analysis yields significant insights:
\begin{itemize}
\item Paradigm of modeling, which can restrain the capability of a model in capturing information contained in the data that reflects physical meanings of reality, has to match the advancements in observation, computation and communication that have dramatically increased the information capacity and diversity of data in describing reality.
\item The dramatically increasing big data has pushed people to a stage to accept and model phenomenon as it exists. It exists in undivided time course. Laws of nature exist in undivided time course~(in a sense of macrophysics). Snapshot model, or model independent of time, stems from the intuition w.r.t. how we perceive. Its information capacity is extremely limited, though it fitted the capability in detecting the world in the early age.
\end{itemize}

Models of load, supply cost and pricing with specific physical meaning and the computational model of pricing have been derived from the general models along the direction of the new paradigm of modeling, named \emph{course-based modeling}.

More interestingly, the general model, though with infinite dimensions, suggests that it is sufficiency rather than completeness that
matters in capturing information and modeling. Actually, the general model should be interpreted as a model system covering models of all dimensions between 1 to $\infty$.
The dimension is a relative concept determined by both the time length chosen to scope a phenomenon and the time length of detecting cycle. An one-hour detecting cycle in an one-day scope yields 24 dimensions. We can adjust the accuracy of the model by adjusting both of the time lengths. By feedback on the adjustment, we are able to evaluate the complexity of a phenomenon in the intended time scope, which contributes to a more effective and efficient modeling solution.

%
%% conference papers do not normally have an appendix
%
%
%% use section* for acknowledgement
%\section*{Acknowledgment}



%\bibliographystyle{IEEEtran}

\end{document}